\documentclass[11pt,a4paper]{article}
\usepackage{amsmath,amssymb}
\usepackage{latexsym}
\usepackage{mathrsfs}
\usepackage{amsfonts}
\usepackage{hyperref}
\usepackage{amsthm}

\hfuzz=\maxdimen
\theoremstyle{definition}
\newtheorem{lemma}{Lemma}[section]
\newtheorem{definition}[lemma]{Definition}
\newtheorem{theorem}[lemma]{Theorem}
\newtheorem{corollary}[lemma]{Corollary}
\newtheorem{remark}{Remark}

\DeclareFixedFont{\Acknowledgment}{OT1}{cmr}{bx}{n}{14pt}
\begin{document}

\title{\bf A new proof for global rigidity of vertex scaling on polyhedral surfaces}
\author{Xu Xu, Chao Zheng}
\maketitle

\begin{abstract}
The vertex scaling for piecewise linear metrics on polyhedral surfaces was introduced by Luo \cite{L1}, who  proved the local rigidity by establishing a variational principle and conjectured the global rigidity.
Luo's conjecture was solved by Bobenko-Pinkall-Springborn \cite{BPS}, who also introduced the vertex scaling for piecewise hyperbolic metrics and proved its global rigidity.
Bobenko-Pinkall-Spingborn's proof is based on their observation of the connection of vertex scaling and the geometry of polyhedra in $3$-dimensional hyperbolic space and the concavity of the volume of ideal and hyper-ideal tetrahedra.
In this paper, we give an elementary and short variational proof
of the global rigidity of vertex scaling without involving $3$-dimensional hyperbolic geometry.
The method is based on continuity of eigenvalues of matrices and the extension of convex functions.
\end{abstract}


\section{Introduction}
The most important two discrete conformal metrics on polyhedral surfaces are circle packing metrics and vertex scaling of polyhedral metrics on surfaces.
There are lots of important works on circle packing metrics, please refer to \cite{An1,An2,BS,CL,DV,K1,L3,MR,S,T1,Z} and others.
In this paper, we focus on vertex scaling of polyhedral metrics on surfaces, which is an analogue of the conformal transformation in Riemannian geometry.
The vertex scaling of piecewise linear metrics (PL metrics for short in the following) on polyhedral surfaces
was introduced physically by R\v{o}cek-Williams \cite{RW} and mathematically by Luo \cite{L1} independently.
Luo proved the local rigidity of vertex scaling for PL metrics by establishing a variational principle
and conjectured the global rigidity in \cite{L1}, where Luo also introduced the corresponding combinatorial Yamabe flow and studied its properties.
Luo's conjecture was solved affirmatively by Bobenko-Pinkall-Springborn in their important work \cite{BPS} by establishing the connection of vertex scaling and the geometry of ideal tetrahedra in $3$-dimensional hyperbolic space and using Rivin's result on the concavity of the volume of ideal tetrahedra \cite{R}.
Bobenko-Pinkall-Springborn \cite{BPS} further introduced the vertex scaling for piecewise hyperbolic metrics (PH metrics for short in the following) and proved its global rigidity by connecting the hyperbolic vertex scaling
to the geometry of hyper-ideal tetrahedra in $3$-dimensional hyperbolic space and using Leibon's result on the concavity of the volume of hyper-ideal tetrahedra \cite{Lei}.
Based on Bobenko-Pinkall-Springborn's observations, the important discrete uniformization theorems for vertex scaling on closed surfaces were recently established in \cite{GGLSW, GLSW, Sp}.
Other related work on the vertex scaling could be found in \cite{GLW, LW, SWGL, Wu, WGS, WZ, X4, XZ1, XZ2,  ZX}.
This paper aims at giving an elementary, direct and short variational proof
for the global rigidity of vertex scaling of PL and PH metrics on surfaces without involving $3$-dimensional hyperbolic geometry.

Suppose $M$ is a closed surface with a triangulation $\mathcal{T}={(V,E,F)}$, where $V,E,F$ represent the sets of vertices, edges and faces respectively.
A discrete metric is a map $l: E\rightarrow (0, +\infty)$ such that the triangle inequalities are satisfied for $l_{ij}, l_{ik}, l_{jk}$ on any triangle $\triangle v_iv_jv_k\in F$,
where $l_{ij}:=l(v_iv_j)$ with $v_iv_j\in E$.
In this case, we can attach a Euclidean metric to each triangle $\triangle v_iv_jv_k\in F$,
which gives rise to a Euclidean triangle, still denoted by $\triangle v_iv_jv_k$.
By gluing the Euclidean triangles isometrically along the edges, we have a PL metric on the triangulated surface $(M, \mathcal{T})$.
If we replace the Euclidean metric by hyperbolic metric, then we obtain a PH metric on the triangulated surface $(M, \mathcal{T})$.
For PL and PH metrics on $(M, \mathcal{T})$, there may exists cone singularities at the vertices, which could be described by combinatorial curvature.
The combinatorial curvature $K_i$ at the vertex $v_i$ is $2\pi$ less the summation of inner angles of triangles at $v_i$.

\begin{definition}[\cite{BPS,L1,RW}]
Suppose $(M, \mathcal{T})$ is a triangulated surface and $u: V\rightarrow \mathbb{R}$ is a function defined on the vertices.
\begin{description}
  \item[(1)] If $l: E\rightarrow (0, +\infty)$ and $\widetilde{l}: E\rightarrow (0, +\infty)$ are two PL metrics on $(M, \mathcal{T})$ with
\begin{equation}\label{Euclidean vertex scaling}
l_{ij}=\widetilde{l}_{ij}e^{\frac{u_i+u_j}{2}}
\end{equation}
for any edge $v_iv_j\in E$,
we say $l$ is a Euclidean vertex scaling of $\widetilde{l}$.
  \item[(2)] If $l: E\rightarrow (0, +\infty)$ and $\widetilde{l}: E\rightarrow (0, +\infty)$ are two PH metrics on $(M, \mathcal{T})$ with
\begin{equation}\label{discrete conformal equivalent formula}
\sinh\frac{l_{ij}}{2}=\sinh \frac{\widetilde{l}_{ij}}{2}e^{\frac{u_i+u_j}{2}}
\end{equation}
for any edge $v_iv_j\in E$,
we say $l$ is a hyperbolic vertex scaling of $\widetilde{l}$.
\end{description}
The function $u: V\rightarrow \mathbb{R}$ is called a discrete conformal factor.
\end{definition}
Bobenko-Pinkall-Spingborn proved the following global rigidity of vertex scaling of polyhedral metrics on surfaces in their important work \cite{BPS}.
\begin{theorem}[\cite{BPS}]\label{main theorem}
Suppose $(M, \mathcal{T})$ is a closed triangulated surface. Then the discrete conformal factor is uniquely determined by the discrete curvature (up to a vector $t(1,1,\cdots, 1), t\in \mathbb{R}$, in the PL case).
\end{theorem}

Let us recall Bobenko-Pinkall-Spingborn's strategy to prove Theorem \ref{main theorem}.
In the PL case, they considered the Legendre transform of the volume of ideal tetrahedra in $3$-dimensional hyperbolic space,
which has an explicit formula in dihedral angles obtained by Milnor \cite{M}.
Based on Rivin's result \cite{R} that the volume of ideal tetrahedra is a concave function of the dihedral angles and could be extended, they extended the definition of Legendre transform of the volume to to be a globally defined convex function.
By modifying the Legendre transform of the volume by a linear function, they showed that this is a globally defined convex extension of Luo's action function which is locally convex.
Then the global rigidity follows from the convexity of the extended function.
In the PH case, the global rigidity is proved similarly with the volume of ideal tetrahedra replaced by
the volume of hyper-ideal tetrahedra with one hyper-ideal vertex and three ideal vertices.
The explicit formula of such hyper-ideal tetrahedra in terms of dihedral angles was obtained by Leibon in \cite{Lei}, where the concavity of the volume was also proved.
Bobenko-Pinkall-Spingborn's approach established the connection of vertex scaling on polyhedral surfaces
and the geometry of hyperbolic polyhedra in $3$-dimensional hyperbolic space.
In this approach they could define the vertex scaling for hyperbolic and spherical polyhedral metrics on surfaces and give the explicit formula of the action functional introduced by Luo, which has lots of applications.

In this paper, we give an elementary, direct and short variational
proof for the global rigidity of vertex scaling of PL and PH metrics on polyhedral surfaces, which does not involve the volume of ideal and hyper-ideal tetrahedra in $3$-dimensional hyperbolic space and their concavity with respect to dihedral angles.
The main idea comes from \cite{BPS, DV, L3,X2}.
The first step is to give a characterization of the admissible space of conformal factors for any given initial discrete metric on a single triangle, which is proved to be simply connected with analytic boundaries by solving a quadratic inequality.
The second step is to prove the Jacobian of the combinatorial curvature with respect to the discrete conformal factor is symmetric and positive definite, which could be reduced to the case that the Jacobian of the inner angles with respect to the conformal factors in a triangle is symmetric and negative definite.
The symmetry could be proved by direct calculations. For the negative definiteness, we introduce a parameterized admissible space of conformal factors, which is the union of admissible spaces of conformal factors on a single triangle with different initial discrete metrics.
This space is proved to be connected, from which the negativity of the Jacobian of inner angles with respect to the conformal factors in a triangle follows easily by the continuity of eigenvalues and calculating at a good point in the parameterized admissible space.
The first step and second step enable us to define a locally convex function on the admissible space of conformal factors for a triangle with fixed initial metric.
The third step is to extend the locally convex function to be a globally defined convex function, from which the global rigidity follows.
This step is accomplished using Luo's extension theorem for continuous closed $1$-forms \cite{L3}, which is a development of Bobenko-Pinkall-Spingborn's extension.
We will give the details of the proof of global rigidity for hyperbolic vertex scaling and just sketch the proof for Euclidean vertex scaling.

The paper is organized as follows.
In Section \ref{Section 2}, we characterize the admissible space of conformal factors for PH metrics on a single triangle.
In Section \ref{Section 3}, we prove the Jacobian matrix of the inner angles in terms of hyperbolic discrete conformal factors in a triangle is symmetric and negative definite, which enables us to define a locally convex function on the admissible space of conformal factors.
In Section \ref{Section 4}, we extend the locally convex function to be a globally defined convex function,
from which the global rigidity of hyperbolic vertex scaling follows.
In Section \ref{Section 5}, we sketch the proof for global rigidity of vertex scaling for PL metrics.

\section{Admissible space of discrete conformal factors for discrete hyperbolic metrics on a triangle}\label{Section 2}

Suppose $\triangle v_iv_jv_k\in F$ is a triangle and $\widetilde{l}_{ij}, \widetilde{l}_{ik}, \widetilde{l}_{jk}$ is a discrete hyperbolic metric on $\triangle v_iv_jv_k$.
The admissible space $\Omega^H_{ijk}(\widetilde{l})$ of discrete conformal factors for the triangle $\triangle v_iv_jv_k$ with discrete hyperbolic metric $\widetilde{l}_{ij}, \widetilde{l}_{ik}, \widetilde{l}_{jk}$
is defined to be the set of discrete conformal factors $({u_i,u_j,u_k})\in \mathbb{R}^3$ such that
the triangle with edge lengths given by formula (\ref{discrete conformal equivalent formula}) exists in 2-dimensional hyperbolic space $\mathbb{H}^2$, i.e.
\begin{equation*}\label{admissible space}
\Omega^{H}_{ijk}(\widetilde{l})=\{(u_i,u_j,u_k)\in \mathbb{R}^3|l_i+l_j>l_k, l_i+l_k>l_j, l_j+l_k>l_i \}.
\end{equation*}
Here and in the following, we use $l_i$ to denote $l_{jk}$ for simplicity.
The parameterized admissible space of conformal factors for the triangle $\triangle v_iv_jv_k$ is defined to be
\begin{equation*}\label{parameterized admissible space}
\Omega^{H}_{ijk}=\{(\widetilde{l}_i, \widetilde{l}_j, \widetilde{l}_k, u_i,u_j,u_k)\in \mathbb{R}^3_{>0}\times \mathbb{R}^3|
l_i+l_j>l_k, l_i+l_k>l_j, l_j+l_k>l_i \},
\end{equation*}
which could be taken as the union of the admissible space $\Omega^{H}_{ijk}(\widetilde{l})$ according
to the parameters given by the initial discrete metrics $(\widetilde{l}_{i}, \widetilde{l}_{j},\widetilde{l}_{k})$.

By formula (\ref{discrete conformal equivalent formula}), if the edge lengths $l_i, l_j, l_k$ satisfy the triangle inequalities, there are some restrictions on the discrete conformal factors.

\begin{lemma}\label{triangle inequalities}
Suppose the triangle $\triangle v_iv_jv_k$ is a triangle with discrete hyperbolic metric $(\widetilde{l}_{i}, \widetilde{l}_{j},\widetilde{l}_{k})$,
$l_i,l_j,l_k$ are the edge lengths defined by formula (\ref{discrete conformal equivalent formula}), then $l_i,l_j,l_k$ satisfy the triangle inequalities if and only if
\begin{equation*}
Q:=-S^4_i\xi^2_i-S^4_j\xi^2_j-S^4_k\xi^2_k+2S^2_iS^2_j\xi_i\xi_j+2S^2_iS^2_k\xi_i\xi_k
+2S^2_jS^2_k\xi_j\xi_k+4S^2_iS^2_jS^2_k>0,
\end{equation*}
where $S_i=\sinh\frac{\widetilde{l}_i}{2},\ \xi_i=e^{-u_i}.$
\end{lemma}
\proof
$l_i,l_j,l_k$ satisfy the triangle inequalities, i.e. $l_i+l_j>l_k, l_i+l_k>l_j, l_j+l_k>l_i,$
which is equivalent to $$\sinh\frac{l_i+l_j+l_k}{2}\sinh\frac{l_i+l_j-l_k}{2}\sinh\frac{l_i+l_k-l_j}{2}\sinh\frac{l_j+l_k-l_i}{2}>0.$$
By direct calculations, we have
\begin{equation*}
\begin{aligned}
&\sinh\frac{l_i+l_j+l_k}{2}\sinh\frac{l_i+l_j-l_k}{2}\sinh\frac{l_i+l_k-l_j}{2}\sinh\frac{l_j+l_k-l_i}{2}\\
=&\frac{1}{4}(\cosh (l_i+l_j)-\cosh l_k)(\cosh l_k-\cosh (l_i-l_j))\\
=&\frac{1}{4}(-\cosh^2 l_i-\cosh^2 l_j-\cosh^2 l_k+2\cosh l_i\cosh l_j\cosh l_k+1)\\
=&-\sinh^4\frac{l_i}{2}-\sinh^4\frac{l_j}{2}-\sinh^4\frac{l_k}{2}+2\sinh^2\frac{l_i}{2}\sinh^2\frac{l_j}{2}\\
&+2\sinh^2\frac{l_i}{2}\sinh^2\frac{l_k}{2}+2\sinh^2\frac{l_j}{2}\sinh^2\frac{l_k}{2}+4\sinh^2\frac{l_i}{2}\sinh^2\frac{l_j}{2}\sinh^2\frac{l_k}{2}\\
=&\xi^{-2}_i\xi^{-2}_j\xi^{-2}_k(-S^4_i\xi^2_i-S^4_j\xi^2_j-S^4_k\xi^2_k+2S^2_iS^2_j\xi_i\xi_j+2S^2_iS^2_k\xi_i\xi_k
+2S^2_jS^2_k\xi_j\xi_k+4S^2_iS^2_jS^2_k),
\end{aligned}
\end{equation*}
where the formula (\ref{discrete conformal equivalent formula}) is used in the last equality.
\qed

Set
\begin{equation}\label{hi hj hk}
\begin{aligned}
h_i&=-S^4_i\xi_i+S^2_iS^2_j\xi_j+S^2_iS^2_k\xi_k,\\
h_j&=-S^4_j\xi_j+S^2_iS^2_j\xi_i+S^2_jS^2_k\xi_k,\\
h_k&=-S^4_k\xi_k+S^2_iS^2_k\xi_i+S^2_jS^2_k\xi_j,
\end{aligned}
\end{equation}
then we have
$$Q=\xi_ih_i+\xi_jh_j+\xi_kh_k+4S^2_iS^2_jS^2_k.$$
Lemma \ref{triangle inequalities} implies that
$(u_i,u_j,u_k)\in \mathbb{R}^3$ is a degenerate hyperbolic discrete conformal factor for a triangle $\triangle v_iv_jv_k$ if and only if $$Q=\xi_ih_i+\xi_jh_j+\xi_kh_k+4S^2_iS^2_jS^2_k\leq{0}.$$
Note that $4S^2_iS^2_jS^2_k>0$ for any $(\widetilde{l}_{i}, \widetilde{l}_{j},\widetilde{l}_{k})\in \mathbb{R}^3_{>0}.$
If $(u_i,u_j,u_k)\in \mathbb{R}^3$ is a degenerate hyperbolic discrete conformal factor, we have $$\xi_ih_i+\xi_jh_j+\xi_kh_k<0,$$
which implies that at least one of $h_i, h_j, h_k$ is negative.
We further have the following result on the signs of $h_i, h_j$ and $h_k$.

\begin{lemma}\label{one negative two positive}
Suppose $(u_i,u_j,u_k)\in \mathbb{R}^3$ is a degenerate hyperbolic discrete conformal factor for a triangle $\triangle v_iv_jv_k$,
then one of $h_i,h_j,h_k$ is negative and the others are positive.
\end{lemma}
\proof
We claim that there exists no subset $\{r,s\}\subset \{i,j,k\}$
such that $h_r\leq 0$ and $h_s\leq 0$, from which the conclusion of the lemma follows.
Otherwise, without loss of generality, we assume $h_i\leq 0,h_j\leq 0$,
which is equivalent to
$S^2_i\xi_i\geq S^2_j\xi_j+S^2_k\xi_k, S^2_j\xi_j\geq S^2_i\xi_i+S^2_k\xi_k.$
This is impossible.
\qed

There is a nice geometric explanation of the result in Lemma \ref{one negative two positive}
in the Euclidean case in terms of circumcircle center. Please refer to Remark \ref{Euclidean geometric explanation}.

\begin{theorem}\label{simple connectedness}
Given any initial nondegenerate hyperbolic discrete metric
$\widetilde{l}=(\widetilde{l}_i, \widetilde{l}_j, \widetilde{l}_k)$ on a triangle $\triangle v_iv_jv_k$,
the admissible space $\Omega^H_{ijk}(\widetilde{l})$ of hyperbolic discrete conformal factors $(u_i,u_j,u_k)\in \mathbb{R}^3$ for the triangle $\triangle v_iv_jv_k$ is nonempty and simply connected.
Furthermore, the set of degenerate hyperbolic discrete conformal factors is
a disjoint union $\bigcup_{\alpha\in \Lambda}{V_\alpha}$, where $\Lambda=\{i,j,k\}$
and $V_\alpha$ is a closed region in $\mathbb{R}^3$ bounded by an analytic graph on $\mathbb{R}^2$.
\end{theorem}
\proof
Suppose $(u_i,u_j,u_k)\in \mathbb{R}^3$ is a degenerate hyperbolic discrete conformal factor for the triangle $\triangle v_iv_jv_k$,
which is equivalent to $Q\leq 0$.
Then by Lemma \ref{one negative two positive}, one of $h_i,h_j,h_k$ is negative and the others are positive.
Without loss of generality, we assume $h_i<0,\ h_j>0,\ h_k>0$.
Note that $Q\leq{0}$ is equivalent to the following quadratic inequality of $\xi_i$
\begin{equation}\label{quadratic inequality}
A_i\xi^2_i+B_i\xi_i+C_i\geq{0},
\end{equation}
where
\begin{equation}\label{coefficient of quadratic inequality}
\begin{split}
A_i&=S^4_i>0,\\
B_i&=-2S^2_i(S^2_j\xi_j+S^2_k\xi_k)<0,\\
C_i&=S^4_j\xi^2_j+S^4_k\xi^2_k-2S^2_jS^2_k\xi_j\xi_k-4S^2_iS^2_jS^2_k.
\end{split}
\end{equation}
By direct calculations, $\Delta_i=B^2_i-4A_iC_i$ is given by
\begin{equation}\label{Delta}
\Delta_i=16S^4_iS^2_jS^2_k\xi_j\xi_k+16S^6_iS^2_jS^2_k>0.
\end{equation}
Combining formula (\ref{quadratic inequality}), (\ref{coefficient of quadratic inequality}) with (\ref{Delta}), we have
\begin{equation*}
\xi_i\geq\frac{-B_i+\sqrt{\Delta_i}}{2A_i} \quad \text{or} \quad\xi_i\leq\frac{-B_i-\sqrt{\Delta_i}}{2A_i}.
\end{equation*}
Note that $-2h_i=2A_i\xi_i+B_i$, so $h_i<0$ is equivalent to $\xi_i>\frac{-B_i}{2A_i}$, which implies $\xi_i\geq\frac{-B_i+\sqrt{\Delta_i}}{2A_i}$. Therefore, $\mathbb{R}^3\setminus \Omega^H_{ijk}(\widetilde{l})\subseteq \bigcup_{\alpha\in \Lambda}{V_\alpha}$,
where
\begin{equation*}
V_i=\left\{(u_i, u_j, u_k)\in \mathbb{R}^3|\xi_i\geq\frac{-B_i+\sqrt{\Delta_i}}{2A_i}\right\}
\end{equation*}
and $V_j, V_k$ are defined similarly.

On the other hand, for any $(u_i,u_j,u_k)\in V_i$, we have $A_i\xi^2_i+B_i\xi_i+C_i\geq{0}$, which implies $Q\leq 0$,
thus $V_i\subseteq \mathbb{R}^3\setminus \Omega^H_{ijk}(\widetilde{l})$. Similar arguments imply $V_j, V_k\subseteq \mathbb{R}^3\setminus \Omega^H_{ijk}(\widetilde{l})$.
Therefore, $\Omega^H_{ijk}(\widetilde{l})=\mathbb{R}^3\backslash\bigcup_{\alpha\in \Lambda}{V_\alpha}$, where $\Lambda=\{i,j,k\}.$

For any $(u_i,u_j,u_k)\in V_i$, we have $\xi_i>\frac{-B_i}{2A_i}$, which is equivalent to $h_i<0$.
Similarly, for $(u_i,u_j,u_k)\in V_j$, we have $h_j<0$ and for $(u_i,u_j,u_k)\in V_k$, we have $h_k<0$.
Then Lemma \ref{one negative two positive} implies $V_i\bigcap V_j=\varnothing$, $V_i\bigcap V_k=\varnothing$, $V_j\bigcap V_k=\varnothing$.

Note that $V_i$ is bounded by an analytic graph on $\mathbb{R}^2$. In fact,
$$V_i=\{(u_i,u_j,u_k)\in \mathbb{R}^3|u_i\leq\log{\frac{2A_i}{-B_i+\sqrt\Delta_i}\}}.$$
This implies $\Omega^H_{ijk}(\widetilde{l})=\mathbb{R}^3\backslash\bigcup_{\alpha\in \Lambda}{V_\alpha}$ is homotopy equivalent to $\mathbb{R}^3$. Therefore, $\Omega^H_{ijk}(\widetilde{l})$ is simply connected.
\qed

\begin{remark}
The ideal of the proof of Theorem \ref{simple connectedness} comes from \cite{X0},
where the first author introduced the method of homotopy deformation to prove the admissible space of sphere packing metrics for a single tetrahedron is simply connected.
This method is then developed \cite{X2} to prove the admissible space of inversive distance packing metrics for a single triangle is simply connected and used to proved the admissible space of Thurston's sphere packing metrics on a tetrahedron is simply connected \cite{HX1,HX2}.
This method has some other applications in characterizing admissible spaces of discrete metrics, see \cite{X3}.
\end{remark}

Note that $Q$ is a continuous function of $(\widetilde{l}_i, \widetilde{l}_j, \widetilde{l}_k, u_i,u_j,u_k)\in \mathbb{R}^3_{>0}\times \mathbb{R}^3$
and the space of hyperbolic discrete metrics $(\widetilde{l}_i, \widetilde{l}_j, \widetilde{l}_k)$ satisfying the triangle inequalities is connected.
As a direct corollary of Lemma \ref{triangle inequalities} and Theorem \ref{simple connectedness}, we have the following result on the parameterized
admissible space $\Omega^H_{ijk}$.

\begin{corollary}\label{connectness of parameterized admissible space}
Suppose $\triangle v_iv_jv_k\in F$. Then the parameterized admissible space $\Omega^H_{ijk}$ is connected.
\end{corollary}

Denote $\alpha_i,\alpha_j,\alpha_k$ as the inner angles in the triangle $\triangle v_iv_jv_k$
so that $\alpha_i$ is opposite to the edge of length $l_i$. We further have the following property of
the inner angles on the admissible space $\Omega^H_{ijk}(\widetilde{l})$.

\begin{lemma}\label{constants to be continuous}
The inner angles $\alpha_i,\alpha_j,\alpha_k$ defined for $(u_i,u_j,u_k)\in\Omega^{H}_{ijk}(\widetilde{l})$
could be extended to be continuous functions $\widetilde{\alpha_i},\widetilde{\alpha}_j,\widetilde{\alpha}_k$ defined on $\mathbb{R}^3$.
\end{lemma}
\proof
By Theorem \ref{simple connectedness},
$\Omega^H_{ijk}(\widetilde{l})=\mathbb{R}^3\backslash\bigcup_{\alpha\in \Lambda}{V_\alpha}$, where $\Lambda=\{i,j,k\}$ and
$V_i=\{(u_i,u_j,u_k)\in\mathbb{R}^3|\xi_i\geq\frac{-B_i+\sqrt\triangle_i}{2A_i}\}$.
Then $\partial V_i=\{(u_i,u_j,u_k)\in \mathbb{R}^3|\xi_i=\frac{-B_i+\sqrt\triangle_i}{2A_i}\}$.
Suppose $(u_i,u_j,u_k)\in\Omega^{H}_{ijk}(\widetilde{l})$ tends a point $(\overline{u}_i,\overline{u}_j,\overline{u}_k)\in\partial V_i$.
By the proof of Lemma \ref{triangle inequalities}, we have
\begin{equation*}
\begin{aligned}
4\xi^{-2}_i\xi^{-2}_j\xi^{-2}_kQ=&4\sinh\frac{l_i+l_j+l_k}{2}\sinh\frac{l_i+l_j-l_k}{2}\sinh\frac{l_i+l_k-l_j}{2}\sinh\frac{l_j+l_k-l_i}{2}\\
=&(\cosh (l_i+l_j)-\cosh l_k)(\cosh l_k-\cosh (l_i-l_j))\\
=&\sinh^2 l_i\sinh^2 l_j-(\cosh l_i\cosh l_j-\cosh l_k)^2\\
=&\sinh^2 l_i\sinh^2 l_j-\sinh^2 l_i\sinh^2 l_j\cos^2\alpha_k\\
=&\sinh^2 l_i\sinh^2 l_j\sin^2\alpha_k.
\end{aligned}
\end{equation*}
As $(u_i,u_j,u_k)\in\Omega^{H}_{ijk}(\widetilde{l})$ tends to $(\overline{u}_i,\overline{u}_j,\overline{u}_k)\in\partial V_i$,
we have $Q\rightarrow{0}$, which implies $\alpha_k\rightarrow 0\ \text{or}\ \pi.$
Similarly, we have $\alpha_i, \alpha_j\rightarrow 0$ or $\pi$.

By formula (\ref{partial alpha i partial u j}), we have
\begin{equation}\label{partial alpha j partial u i}
\begin{aligned}
\frac{\partial \alpha_j}{\partial u_i}
=&\frac{\cosh l_i+\cosh l_j-\cosh l_k-1}{A(\cosh l_k+1)}\\
=&\frac{\sinh^2\frac{l_i}{2}+\sinh^2\frac{l_j}{2}-\sinh^2\frac{l_k}{2}}{A(\sinh^2\frac{l_k}{2}+1)}\\
=&\frac{\xi^{-1}_i\xi^{-1}_j\xi^{-1}_k}{A(S^2_k\xi^{-1}_i\xi^{-1}_j+1)}(S^2_i\xi_i+S^2_j\xi_j-S^2_k\xi_k)\\
=&\frac{\xi^{-1}_i\xi^{-1}_j\xi^{-1}_kh_k}{AS^2_k(S^2_k\xi^{-1}_i\xi^{-1}_j+1)},
\end{aligned}
\end{equation}
where $A=\sinh l_j\sinh l_k\sin\alpha_i,\ S_i=\sinh\frac{\widetilde{l}_i}{2},\ \xi_i=e^{-u_i}$ and $h_k$ is defined by formula (\ref{hi hj hk}).
Note that for $(\overline{u}_i,\overline{u}_j,\overline{u}_k)\in\partial V_i$,
by Lemma \ref{one negative two positive} and the proof of Theorem \ref{simple connectedness},
we have $h_i<0,\ h_j>0,\ h_k>0$ at $(\overline{u}_i,\overline{u}_j,\overline{u}_k)$.
By formula (\ref{partial alpha j partial u i}), we have $\frac{\partial\alpha_j}{\partial u_i}>0$ for $(u_i,u_j,u_k)\in\Omega^{H}_{ijk}(\widetilde{l})$ around $(\overline{u}_i,\overline{u}_j,\overline{u}_k)\in\partial V_i$.
This implies $\alpha_j\rightarrow 0$ as $(u_i,u_j,u_k)\rightarrow(\overline{u}_i,\overline{u}_j,\overline{u}_k)\in\partial V_i$. Otherwise, we have $\alpha_j\rightarrow \pi$ as $(u_i,u_j,u_k)\rightarrow(\overline{u}_i,\overline{u}_j,\overline{u}_k)\in\partial V_i$
and then $\frac{\partial\alpha_j}{\partial u_i}>0$ implies $\alpha_j>\pi$ for
some $(u_i,u_j,u_k)\in\Omega^{H}_{ijk}(\widetilde{l})$ around $(\overline{u}_i,\overline{u}_j,\overline{u}_k)\in\partial V_i$, which
is impossible for hyperbolic triangles.
Similarly, we have $\alpha_k\rightarrow 0$ as $(u_i,u_j,u_k)\rightarrow(\overline{u}_i,\overline{u}_j,\overline{u}_k)\in\partial V_i$.

Furthermore, we have the following formula for the area $S$ of the hyperbolic triangle in terms of the edge lengths (\cite{V} page 66)
\begin{equation*}
\begin{aligned}
\tan^2\frac{S}{4}
=&\tanh\frac{p}{2}\tanh\frac{p-l_i}{2}\tanh\frac{p-l_j}{2}\tanh\frac{p-l_k}{2}\\
=&\frac{\xi^{-2}_i\xi^{-2}_j\xi^{-2}_kQ}{64\cosh^2\frac{p}{2}\cosh^2\frac{p-l_i}{2}\cosh^2\frac{p-l_j}{2}\cosh^2\frac{p-l_k}{2}},
\end{aligned}
\end{equation*}
where $p=\frac{1}{2}(l_i+l_j+l_k)$.
Note that $Q\rightarrow 0$ as $(u_i,u_j,u_k)\rightarrow (\overline{u}_i,\overline{u}_j,\overline{u}_k)\in\partial V_i$,
we have $S\rightarrow 0$.
Then we have $\alpha_i\rightarrow \pi$ as $(u_i,u_j,u_k)\rightarrow (\overline{u}_i,\overline{u}_j,\overline{u}_k)\in\partial V_i$ by $S=\pi-\alpha_i-\alpha_j-\alpha_k$ and
$\alpha_j,\alpha_k\rightarrow 0$.
The case for the boundary $\partial V_j$ and $\partial V_k$ could be discussed similarly.

Therefore, we can extend $\alpha_i,\alpha_j,\alpha_k$ defined on $\Omega^{H}_{ijk}(\widetilde{l})$ to be continuous functions defined on $\mathbb{R}^3$ by setting
\begin{eqnarray*}
\widetilde{\alpha}_i(u_i,u_j,u_k)=
\begin{cases}
\alpha_i,  &\text{if} \ (u_i,u_j,u_k)\in \Omega^{H}_{ijk},\\
\pi,  &\text{if} \ (u_i,u_j,u_k)\in V_i,\\
0,  &\text{if}\ (u_i,u_j,u_k)\in V_j\ \text{or}\ V_k.
\end{cases}
\end{eqnarray*}
This completes the proof of the lemma.
\qed
\begin{remark}
By the proof of Lemma \ref{constants to be continuous}, we have $\frac{\partial\alpha_j}{\partial u_i}\rightarrow +\infty$
and $\frac{\partial\alpha_k}{\partial u_i}\rightarrow +\infty$ as $(u_i,u_j,u_k)\rightarrow (\overline{u}_i,\overline{u}_j,\overline{u}_k)\in\partial V_i$.
Recall the following formula obtained by Glickenstein-Thomas (\cite{GT} Proposition 9)
\begin{equation*}
\frac{\partial S}{\partial u_i}=\frac{\partial \alpha_j}{\partial u_i}(\cosh l_{k}-1)+\frac{\partial \alpha_k}{\partial u_i}(\cosh l_{j}-1),
\end{equation*}
where $S$ is the area of $\triangle v_iv_jv_k$, which could also be proved using Lemma \ref{the elements of matrix} directly.
For hyperbolic vertex scaling,
we have $\frac{\partial S}{\partial u_i}\rightarrow +\infty$, which implies
\begin{equation*}
\frac{\partial \alpha_i}{\partial u_i}=-\frac{\partial S}{\partial u_i}-\frac{\partial \alpha_j}{\partial u_i}-\frac{\partial \alpha_k}{\partial u_i}
\rightarrow -\infty
\end{equation*}
as $(u_i,u_j,u_k)\rightarrow (\overline{u}_i,\overline{u}_j,\overline{u}_k)\in\partial V_i$.
\end{remark}

\section{Negative definiteness of Jacobian matrix}\label{Section 3}

\begin{lemma}\label{the elements of matrix}
For any triangle $\triangle v_iv_jv_k$, let $l_i,l_j,l_k$ be edge lengths of a hyperbolic triangle and $\alpha_i,\alpha_j,\alpha_k$ be the opposite angles so that $\alpha_i$ is facing the edge of length $l_i$, then
\begin{equation}\label{partial alpha i partial u j}
\frac{\partial \alpha_i}{\partial u_j}=\frac{\partial \alpha_j}{\partial u_i}=\frac{\cosh l_i+\cosh l_j-\cosh l_k-1}{A(\cosh l_k+1)},\\
\end{equation}
\begin{equation*}
\begin{split}
\frac{\partial \alpha_i}{\partial u_i}
=&\frac{\cosh^2 l_j+\cosh^2 l_k-2\cosh l_i\cosh l_j\cosh l_k+(1-\cosh l_i)(\cosh l_j+\cosh l_k)}{A(1+\cosh l_j)(1+\cosh l_k)},
\end{split}
\end{equation*}
where $A=\sinh l_j\sinh l_k\sin\alpha_i$.
\end{lemma}
\proof
By the derivative cosine law (see Lemma A1 in \cite{CL} for example), we have
\begin{equation*}
\begin{aligned}
\frac{\partial \alpha_i}{\partial l_i}=\frac{\sinh l_i}{A},\
\frac{\partial \alpha_i}{\partial l_j}=\frac{-\sinh l_i\cos\alpha_k}{A},\
\frac{\partial \alpha_i}{\partial l_k}=\frac{-\sinh l_i\cos\alpha_j}{A},
\end{aligned}
\end{equation*}
where $A=\sinh l_j\sinh l_k\sin\alpha_i$.
By formula (\ref{discrete conformal equivalent formula}), we have
$$\frac{\partial l_i}{\partial u_i}=0,\ \frac{\partial l_i}{\partial u_j}=\frac{\partial l_i}{\partial u_k}=\tanh\frac{l_i}{2}.$$
Then according to the chain rules, we have
\begin{equation*}
\begin{aligned}
\frac{\partial \alpha_i}{\partial u_j}
=&\frac{\partial \alpha_i}{\partial l_i}\frac{\partial l_i}{\partial u_j}+\frac{\partial \alpha_i}{\partial l_j}\frac{\partial l_j}{\partial u_j}+\frac{\partial \alpha_i}{\partial l_k}\frac{\partial l_k}{\partial u_j}\\
=&\frac{\sinh l_i}{A}\tanh\frac{l_i}{2}-\frac{\sinh l_i\cos\alpha_j}{A}\tanh\frac{l_k}{2}\\
=&\frac{\sinh l_i}{A}\frac{\sinh l_i}{1+\cosh l_i}-\frac{\sinh l_i\cos\alpha_j}{A}\frac{\sinh l_k}{1+\cosh l_k}\\
=&\frac{\cosh l_i+\cosh l_j-\cosh l_k-1}{A(\cosh l_k+1)},
\end{aligned}
\end{equation*}
which implies $\frac{\partial \alpha_i}{\partial u_j}=\frac{\partial \alpha_j}{\partial u_i}$.
Similarly, we have
\begin{equation*}
\begin{aligned}
\frac{\partial \alpha_i}{\partial u_i}
=&\frac{\partial \alpha_i}{\partial l_i}\frac{\partial l_i}{\partial u_i}+\frac{\partial \alpha_i}{\partial l_j}\frac{\partial l_j}{\partial u_i}+\frac{\partial \alpha_i}{\partial l_k}\frac{\partial l_k}{\partial u_i}\\
=&\frac{\cosh^2 l_j+\cosh^2 l_k-2\cosh l_i\cosh l_j\cosh l_k+(1-\cosh l_i)(\cosh l_j+\cosh l_k)}{A(1+\cosh l_j)(1+\cosh l_k)}.
\end{aligned}
\end{equation*}
\qed

Lemma \ref{the elements of matrix} shows that the matrix
\begin{equation*}
\begin{aligned}
\Lambda^H_{ijk}=\frac{\partial (\alpha_i, \alpha_j, \alpha_k)}{\partial (u_i, u_j, u_k)}
=\left(
   \begin{array}{ccc}
     \frac{\partial \alpha_i}{\partial u_i} & \frac{\partial \alpha_i}{\partial u_j} & \frac{\partial \alpha_i}{\partial u_k} \\
     \frac{\partial \alpha_j}{\partial u_i} & \frac{\partial \alpha_j}{\partial u_j} & \frac{\partial \alpha_j}{\partial u_k} \\
     \frac{\partial \alpha_k}{\partial u_i} & \frac{\partial \alpha_k}{\partial u_j} & \frac{\partial \alpha_k}{\partial u_k} \\
   \end{array}
 \right)
\end{aligned}
\end{equation*}
is symmetric on $\Omega^H_{ijk}$. Furthermore, one has the following result for the matrix $\Lambda^H_{ijk}$.

\begin{theorem}[\cite{BPS}]\label{negative definite of matrix}
The matrix $\Lambda^H_{ijk}$ is symmetric, negative definite on $\Omega^H_{ijk}$.
\end{theorem}
\proof
By the chain rules, we have
\begin{equation*}
\begin{aligned}
\Lambda^H_{ijk}=\frac{\partial (\alpha_i, \alpha_j, \alpha_k)}{\partial (u_i, u_j, u_k)}
=\frac{\partial (\alpha_i, \alpha_j, \alpha_k)}{\partial (l_i, l_j, l_k)}\cdot \frac{\partial (l_i, l_j, l_k)}{\partial (u_i, u_j, u_k)}.
\end{aligned}
\end{equation*}
By the calculations in the proof of Lemma \ref{the elements of matrix}, we have
\begin{equation}\label{matrix in proof 1}
\begin{aligned}
\frac{\partial (\alpha_i, \alpha_j, \alpha_k)}{\partial (l_i, l_j, l_k)}
=-\frac{1}{A}
 \left(
   \begin{array}{ccc}
     \sinh l_i & 0 & 0 \\
     0 & \sinh l_j & 0 \\
     0 & 0 & \sinh l_k \\
   \end{array}
 \right)
 \left(
   \begin{array}{ccc}
     -1 & \cos\alpha_k & \cos\alpha_j \\
     \cos\alpha_k & -1 & \cos\alpha_i \\
     \cos\alpha_j & \cos\alpha_i & -1 \\
   \end{array}
 \right)\\
\end{aligned}
\end{equation}
and
\begin{equation}\label{matrix in proof 2}
\begin{aligned}
\frac{\partial (l_i, l_j, l_k)}{\partial (u_i, u_j, u_k)}
=  \left(
   \begin{array}{ccc}
     \tanh\frac{l_i}{2} & 0 & 0 \\
     0 & \tanh\frac{l_j}{2} & 0 \\
     0 & 0 & \tanh\frac{l_k}{2} \\
   \end{array}
   \right)
   \left(
     \begin{array}{ccc}
       0 & 1 & 1 \\
       1 & 0 & 1 \\
       1 & 1 & 0 \\
     \end{array}
   \right).
\end{aligned}
\end{equation}
Denote the last matrices in (\ref{matrix in proof 1}) and (\ref{matrix in proof 2}) as $\Phi$ and $\mathcal{R}$ respectively.
By direct calculations, we have
\begin{equation*}
\begin{aligned}
\det \Phi=&-1+\cos\alpha_i^2+\cos\alpha_j^2+\cos\alpha_k^2+2\cos\alpha_i\cos\alpha_j\cos\alpha_k\\
=&4\cos\frac{\alpha_i+\alpha_j-\alpha_k}{2}\cos\frac{\alpha_i-\alpha_j+\alpha_k}{2}
 \cos\frac{\alpha_i+\alpha_j+\alpha_k}{2}\cos\frac{\alpha_i-\alpha_j-\alpha_k}{2}>0,\\
\det \mathcal{R}=&2>0
\end{aligned}
\end{equation*}
for any $(\widetilde{l}_i, \widetilde{l}_j, \widetilde{l}_k, u_i,u_j,u_k)\in \Omega^H_{ijk}$,
which implies $\det \Lambda^H_{ijk}<0$ and then the Jacobian matrix $\Lambda^H_{ijk}$ is non-singular.
Therefore, the eigenvalues of $\Lambda^H_{ijk}$ are non-zero.
Combining with the continuity of the eigenvalues and the connectivity of the parameterized admissible space
$\Omega^H_{ijk}$ in Corollary \ref{connectness of parameterized admissible space},
the eigenvalues of $\Lambda^H_{ijk}$ never change signs. So we just need to calculate at one point in $\Omega^H_{ijk}$ to prove that
the eigenvalues of $\Lambda^H_{ijk}$ are negative and then $\Lambda^H_{ijk}$ is negative definite.
Note that $p=(1, 1, 1, 0, 0, 0)\in \Omega^H_{ijk}$.
By Lemma \ref{the elements of matrix}, we have
\begin{equation*}
\begin{aligned}
\Lambda^H_{ijk}(p)
=&\frac{-(\cosh 1-1)}{A(1+\cosh 1)}\left(
                                        \begin{array}{ccc}
                                          2\cosh 1 & -1 & -1 \\
                                          -1 & 2\cosh 1 & -1 \\
                                          -1 & -1 & 2\cosh 1 \\
                                        \end{array}
                                      \right),
\end{aligned}
\end{equation*}
which is negative definite.
Therefore, the eigenvalues of the Jacobian matrix $\Lambda^H_{ijk}$ at $p=(1, 1, 1, 0, 0, 0)$
are negative.
This completes the proof of the theorem.
\qed

\begin{remark}
Theorem \ref{negative definite of matrix} was first obtained by Bobenko-Pinkall-Springborn in their important work \cite{BPS} by taking the Jacobian matrix $\Lambda^H_{ijk}$ as the Hessian matrix of Legendre transform of
the volume of hyper-ideal tetrahedra in $3$-dimensional hyperbolic space with prescribed metric.
The negativity of $\Lambda^H_{ijk}$ follows from Leibon's concavity of the volume of hyper-ideal tetrahedra with one hyper-ideal
vertex and three ideal vertices \cite{Lei}, which depends on the explicit form of the volume formula in terms of dihedral angles.
The proof of Theorem \ref{negative definite of matrix} presented here involves
only the cosine law and the continuity of the eigenvalues.
\end{remark}

\section{Proof of the global rigidity of hyperbolic vertex scaling}\label{Section 4}

By Theorem \ref{simple connectedness} and Theorem \ref{negative definite of matrix},
the following function
$$F_{ijk}(u_i,u_j,u_k)=\int_{(\overline{u}_i,\overline{u}_j,\overline{u}_k)}^{(u_i,u_j,u_k)}\alpha_idu_i+\alpha_jdu_j+\alpha_kdu_k$$
is a well-defined locally strictly concave function of $(u_i,u_j,u_k)\in \Omega^{H}_{ijk}(\widetilde{l})$.
We need to extend $F_{ijk}$ to be a globally defined concave function on $\mathbb{R}^3$.
Recall the following definition of closed continuous $1$-form and
extension of locally convex function of Luo \cite{L3}, which is a development of Bobenko-Pinkall-Spingborn's extension in \cite{BPS}.
\begin{definition}[\cite{L3}, Definition 2.3]
A differential 1-form $w=\sum_{i=1}^n a_i(x)dx^i$ in an open set $U\subset \mathbb{R}^n$ is said to be continuous if each $a_i(x)$ is continuous on $U$. A continuous differential 1-form $w$ is called closed if $\int_{\partial \tau}w=0$ for each
triangle $\tau\subset U$.
\end{definition}

\begin{theorem}[\cite{L3}, Corollary 2.6]\label{Luo's convex extention}
Suppose $X\subset \mathbb{R}^n$ is an open convex set and $A\subset X$ is an open subset of $X$ bounded by a real analytic codimension-1 submanifold in $X$. If $w=\sum_{i=1}^na_i(x)dx_i$ is a continuous closed 1-form on $A$ so that $F(x)=\int_a^x w$ is locally convex on $A$ and each $a_i$ can be extended continuous to $X$ by constant functions to a function $\widetilde{a}_i$ on $X$, then  $\widetilde{F}(x)=\int_a^x\sum_{i=1}^n\widetilde{a}_i(x)dx_i$ is a $C^1$-smooth
convex function on $X$ extending $F$.
\end{theorem}

By Lemma \ref{constants to be continuous} and Theorem \ref{Luo's convex extention},
$F_{ijk}(u_i,u_j,u_k)$ defined on $\Omega^{H}_{ijk}(\widetilde{l})$ could be extended to be the following function
\begin{equation*}\label{global hyperbolic energy function}
\widetilde{F}_{ijk}(u_i,u_j,u_k)=\int_{(\overline{u}_i,\overline{u}_j,\overline{u}_k)}^
{(u_i,u_j,u_k)}\widetilde{\alpha}_idu_i+\widetilde{\alpha_j}du_j+\widetilde{\alpha}_kdu_k,
\end{equation*}
which is a $C^1$-smooth concave function defined on $\mathbb{R}^3$ with $\nabla_u \widetilde{F}_{ijk}=(\widetilde{\alpha}_i,\widetilde{\alpha_j},\widetilde{\alpha}_k)^{T}.$
Set
\begin{equation*}
\widetilde{F}(u_1,\cdots,u_{|V|})=-\sum_{\triangle v_iv_jv_k\in F}\widetilde{F}_{ijk}(u_i,u_j,u_k)+\int_{\overline{u}}^{u}2\pi\sum_{i=1}^{|V|} du_i,
\end{equation*}
where $|V|$ is the number of vertices. Then $\widetilde{F}(u_1,\cdots,u_{|V|})$ is a $C^1$ smooth convex function on $\mathbb{R}^V$ with
\begin{equation*}\label{nabla u_i}
\nabla_{u_i}\widetilde{F}(u_1,\cdots,u_{|V|})=-\sum_{\triangle ijk\in F}\widetilde \alpha_i+2\pi=\widetilde{K}_i,
\end{equation*}
where $\widetilde{K}_i=2\pi-\sum_{\triangle v_iv_jv_k\in F}\widetilde \alpha_i$ is an extension of $K_i$.
Then the global rigidity of hyperbolic vertex scaling follows from the convexity of $\widetilde{F}$ on $\mathbb{R}^V$
and the locally strict convexity of $\widetilde{F}$ on $\cap_{\triangle v_iv_jv_k\in F}\Omega_{ijk}^H(\widetilde{l})$.
This completes the proof of Theorem \ref{main theorem} in the hyperbolic case.

\section{Rigidity for vertex scaling of PL metrics}\label{Section 5}
As the main steps for the proof of global rigidity of Euclidean vertex scaling is paralleling to the hyperbolic case, we just list the main steps here.

Given any initial discrete Euclidean metric $\widetilde{l}_{ij}, \widetilde{l}_{ik}, \widetilde{l}_{jk}$ on the triangle $\triangle v_iv_jv_k$, the admissible space $\Omega^{E}_{ijk}(\widetilde{l})$ of Euclidean conformal factors is defined to be
\begin{equation*}
\Omega^{E}_{ijk}(\widetilde{l})=\{(u_i,u_j,u_k)\in \mathbb{R}^3|l_i+l_j>l_k, l_i+l_k>l_j, l_j+l_k>l_i \},
\end{equation*}
where the edge lengths are given by formula (\ref{Euclidean vertex scaling}) and
we use $l_i$ to denote $l_{jk}$ for simplicity.
The Euclidean parameterized admissible space of conformal factors for the triangle $\triangle v_iv_jv_k$ is defined to be
\begin{equation*}
\Omega^{E}_{ijk}=\{(\widetilde{l}_i, \widetilde{l}_j, \widetilde{l}_k, u_i,u_j,u_k)\in \mathbb{R}^3_{>0}\times \mathbb{R}^3|l_i+l_j>l_k, l_i+l_k>l_j, l_j+l_k>l_i \}.
\end{equation*}

\begin{lemma}\label{Euclidean triangle inequalities}
Suppose the triangle $\triangle v_iv_jv_k$ is a topological triangle,
$l_i,l_j,l_k$ are the edge lengths defined by (\ref{Euclidean vertex scaling}), then the triangle inequalities are satisfied if and only if
\begin{equation*}
Q:=-\widetilde{l}^4_i\xi^2_i-\widetilde{l}^4_j\xi^2_j-\widetilde{l}^4_k\xi^2_k+2\widetilde{l}^2_i\widetilde{l}^2_j\xi_i\xi_j
+2\widetilde{l}^2_i\widetilde{l}^2_k\xi_i\xi_k+2\widetilde{l}^2_j\widetilde{l}^2_k\xi_j\xi_k> 0,
\end{equation*}
where $\xi_i=e^{-u_i}.$
\end{lemma}
Set
\begin{equation}\label{Euclidean hi hj hk}
\begin{aligned}
h_i&=-\widetilde{l}^4_i\xi_i+\widetilde{l}^2_i\widetilde{l}^2_j\xi_j+\widetilde{l}^2_i\widetilde{l}^2_k\xi_k,\\
h_j&=-\widetilde{l}^4_j\xi_j+\widetilde{l}^2_i\widetilde{l}^2_j\xi_i+\widetilde{l}^2_j\widetilde{l}^2_k\xi_k,\\
h_k&=-\widetilde{l}^4_k\xi_k+\widetilde{l}^2_i\widetilde{l}^2_k\xi_i+\widetilde{l}^2_j\widetilde{l}^2_k\xi_j.\\
\end{aligned}
\end{equation}
Then we have
$Q=\xi_ih_i+\xi_jh_j+\xi_kh_k.$
$(u_i,u_j,u_k)\in \mathbb{R}^3$ is a degenerate Euclidean discrete conformal factor if and only if
$Q=\xi_ih_i+\xi_jh_j+\xi_kh_k\leq 0$, which implies at least one of $h_i, h_j, h_k$ is nonpositive.
Similar to Lemma \ref{one negative two positive} in the hyperbolic case, we have the following result on the signs of $h_i, h_j, h_k$
in the Euclidean case.

\begin{lemma}\label{Euclidean one negative two positive}
Suppose $(u_i,u_j,u_k)\in \mathbb{R}^3$ is a degenerate Euclidean discrete conformal factor for a triangle $\triangle v_iv_jv_k$, then one of $h_i,h_j,h_k$ is negative and the others are positive.
\end{lemma}

\begin{remark}\label{Euclidean geometric explanation}
Lemma \ref{Euclidean one negative two positive} has the following interesting geometrical explanation.
For a Euclidean triangle $\triangle v_iv_jv_k$ with a nondegenerate discrete conformal factor $(u_i,u_j,u_k)$, there exists a geometric center $C_{ijk}$ (\cite{G} Proposition 4) of the triangle $\triangle v_iv_jv_k$ with the same Euclidean distance from $C_{ijk}$ to each vertex of the triangle, which is in fact the circumcircle center for vertex scaling of PL metrics.
$h_i$ in formula (\ref{Euclidean hi hj hk}) is positive multiplication of the signed distance $h_{jk,i}$ from $C_{ijk}$ to the edge $\{jk\}$, which is defined to be positive if $C_{ijk}$ is on the same side of the line determined by $\{jk\}$ as the triangle $\triangle v_iv_jv_k$ and negative otherwise (or zero if $C_{ijk}$ is on the edge). By direct calculations, we have the following relationship for $h_i$ and $h_{jk,i}$
\begin{equation*}
h_{jk,i}=\frac{\xi^{-1}_i\xi^{-\frac{3}{2}}_j\xi^{-\frac{3}{2}}_k}{8S\widetilde{l}_i}h_i,
\end{equation*}
where $S$ is the area of the Euclidean triangle $\triangle v_iv_jv_k$.
See \cite{G3} for more general cases.
For degenerate
conformal factors for Euclidean vertex scaling, Lemma $\ref{Euclidean one negative two positive}$ implies that
the circumcircle center lies in some special regions in the plane relative to the triangle $\triangle v_iv_jv_k$.
\end{remark}

\begin{theorem}[\cite{L1}]\label{Euclidean simple connectedness}
Given any initial nondegenerate Euclidean discrete metric
$\widetilde{l}=(\widetilde{l}_i, \widetilde{l}_j, \widetilde{l}_k)$ on a triangle $\triangle v_iv_jv_k$,
the admissible space $\Omega^E_{ijk}(\widetilde{l})$ of Euclidean discrete conformal factors $(u_i,u_j,u_k)\in \mathbb{R}^3$ for the triangle $\triangle v_iv_jv_k$ is nonempty and simply connected. Furthermore, the set of degenerate Euclidean discrete conformal factors is a disjoint union $\bigcup_{\alpha\in \Lambda}{V_\alpha}$, where $\Lambda=\{i,j,k\}$ and $V_\alpha$ is bounded by an analytic graph on $\mathbb{R}^2$ with
\begin{equation*}
V_i=\{(u_i,u_j,u_k)\in \mathbb{R}^3| u_i\leq -\ln (\widetilde{l}^2_je^{-u_j}+\widetilde{l}^2_ke^{-u_k})+2\ln \widetilde{l}_i \}.
\end{equation*}
As a corollary, $\Omega^{E}_{ijk}$ is connected.
\end{theorem}

Following the hyperbolic case, as an application of Theorem \ref{Euclidean simple connectedness},
we have the following result, which was obtained by Luo \cite{L1} by direct calculations.

\begin{theorem}[\cite{L1}]\label{Euclidean negative definite of matrix}
The matrix $\Lambda^E_{ijk}=[\frac{\partial\alpha_r}{\partial u_s}]_{3\times 3}$ is symmetric, semi-negative definite on $\Omega^E_{ijk}(\widetilde{l})$ with null space $\{(t,t,t)\in \mathbb{R}^3|t\in \mathbb{R}\}$.
\end{theorem}

\begin{remark}
In fact, by the derivative cosine law (see \cite{CL} for example), we have
$\frac{\partial\alpha_i}{\partial l_i}=\frac{l_i}{2S},\ \frac{\partial\alpha_i}{\partial l_j}=-\frac{l_i\cos \alpha_k}{2S},\ \frac{\partial\alpha_i}{\partial l_k}=-\frac{l_i\cos \alpha_j}{2S},$
where $S$ is the area of the Euclidean triangle $\triangle v_iv_jv_k$.
According to formula (\ref{Euclidean vertex scaling}), we have
$\frac{\partial l_i}{\partial u_i}=0,\ \frac{\partial l_i}{\partial u_j}=\frac{\partial l_i}{\partial u_k}=\frac{l_i}{2}.$
By direct calculations with the chain rules, we have
$
\frac{\partial\alpha_i}{\partial u_j}
=\frac{l_il_j\cos\alpha_k}{4S},\ \
\frac{\partial\alpha_i}{\partial u_i}
=-\frac{l^2_i}{4S},
$
which implies $\Lambda^E_{ijk}=[\frac{\partial\alpha_r}{\partial u_s}]_{3\times 3}$ is symmetric.
Luo \cite{L1} proved the semi-negative definiteness of $\Lambda^E_{ijk}$ by direct calculations for any nondegenerate Euclidean conformal factor.
If we use the connectivity of $\Omega_{ijk}^E$, we just need to check the signs of the eigenvalues of $\Lambda^E_{ijk}$
at the point $p=(1,1,1,0,0,0)\in \Omega_{ijk}^E$.
By direct calculations, we have
\begin{equation*}
\begin{aligned}
\Lambda^E_{ijk}(p)=\frac{-\sqrt{3}}{6}\left(
                                        \begin{array}{ccc}
                                          2 & -1 & -1 \\
                                          -1 & 2 & -1 \\
                                          -1 & -1 & 2 \\
                                        \end{array}
                                      \right),
\end{aligned}
\end{equation*}
which has two negative eigenvalues and one zero eigenvalue.
This also implies semi-negative definiteness of $\Lambda^E_{ijk}$.
\end{remark}

\begin{lemma}[\cite{BPS}]\label{Euclidean extension}
Suppose $(u_i,u_j,u_k)\in \mathbb{R}^3$ is a nondegenerate Euclidean discrete conformal factor for a triangle $\triangle v_iv_jv_k$, denote $\alpha_i$ as the angle at vertex $v_i$. Then $\alpha_i,\alpha_j,\alpha_k$ defined for $(u_i,u_j,u_k)\in\Omega^{E}_{ijk}(\widetilde{l})$ could be extended by constants to be continuous functions $\widetilde{\alpha_i},\widetilde{\alpha}_j,\widetilde{\alpha}_k$ defined on $\mathbb{R}^3$.
\end{lemma}

By Theorem \ref{Euclidean simple connectedness} and Theorem \ref{Euclidean negative definite of matrix},
the following function
$$F_{ijk}(u_i,u_j,u_k)=\int_{(\overline{u}_i,\overline{u}_j,\overline{u}_k)}^{(u_i,u_j,u_k)}\alpha_idu_i+\alpha_jdu_j+\alpha_kdu_k$$
is a well-defined locally concave function of $(u_i,u_j,u_k)\in \Omega^{E}_{ijk}(\widetilde{l})$ with $F_{ijk}(u_i+t,u_j+t,u_k+t)=F_{ijk}(u_i,u_j,u_k)+t\pi$.
By Lemma \ref{Euclidean extension} and Theorem \ref{Luo's convex extention},
$F_{ijk}(u_i,u_j,u_k)$ defined on $\Omega^{E}_{ijk}(\widetilde{l})$ could be extended to the following function
\begin{equation*}\label{global Euclidean energy function}
\widetilde{F}_{ijk}(u_i,u_j,u_k)=\int_{(\overline{u}_i,\overline{u}_j,\overline{u}_k)}^
{(u_i,u_j,u_k)}\widetilde{\alpha}_idu_i+\widetilde{\alpha}_jdu_j+\widetilde{\alpha}_kdu_k,
\end{equation*}
which is a $C^1$-smooth concave function defined on $\mathbb{R}^3$ with $\nabla_u \widetilde{F}_{ijk}=(\widetilde{\alpha}_i,\widetilde{\alpha_j},\widetilde{\alpha}_k)^{T}.$
Then the following of the proof for the global rigidity for Euclidean vertex scaling is almost the same
as the hyperbolic case. We omit the details here.

\begin{remark}
In the Euclidean case, similar idea to use Luo's extension theorem \ref{Luo's convex extention} to extend
$F_{ijk}(u_i,u_j,u_k)$ appears in \cite{GJ0}, where the extension depends on the
simply connectivity of the admissible space  $\Omega^{E}_{ijk}(\widetilde{l})$ and
negative semi-definiteness of $\Lambda^E_{ijk}$ obtained by Luo \cite{L1}.
Here we provide a unified approach to prove the simply connectivity of the admissible space of conformal factors
and the negative definiteness of the Jacobian matrix $[\frac{\partial\alpha_r}{\partial u_s}]_{3\times 3}$
for a triangle in the Euclidean and hyperbolic cases.
\end{remark}

\textbf{Acknowledgements}\\[8pt]
The research of the second author is supported by
the Fundamental Research Funds for the Central Universities under
grant no. 2042020kf0199
and National Natural Science Foundation of China under grant no. 61772379.

(Xu Xu) School of Mathematics and Statistics, Wuhan University, Wuhan 430072, P.R. China

E-mail: xuxu2@whu.edu.cn\\[2pt]

(Chao Zheng) School of Mathematics and Statistics, Wuhan University, Wuhan 430072, P.R. China

E-mail: 2019202010023@whu.edu.cn\\[2pt]

\end{document}